\batchmode
\documentclass[final,3p,times]{elsarticle}

\usepackage{lineno,hyperref}
\modulolinenumbers[5]

\journal{}

\usepackage{amsmath}
\usepackage{amsfonts}
\usepackage{amssymb}
\usepackage{dsfont}
\usepackage{graphicx}
\usepackage{epstopdf}
\usepackage{pgfplots}
\usepackage{diagbox}
\usepackage{dirtytalk, float}
\usepackage{multicol}

\usepackage{lipsum}
\makeatletter
\def\ps@pprintTitle{%
 \let\@oddhead\@empty
 \let\@evenhead\@empty
 \def\@oddfoot{}%
 \let\@evenfoot\@oddfoot}
\makeatother

\begin{document}

\begin{frontmatter}

\title{Two-field mixed $hp$-finite elements for time-dependent problems\\in the refined theories of thermodynamics}

\author{Bal\'{a}zs T\'{o}th}
\cortext[mycorrespondingauthor]{Corresponding author}
\address{Institute of Applied Mechanics, University of Miskolc, Miskolc-Egyetemv\'{a}ros, H-3515, Hungary}
\ead{mechtb@uni-miskolc.hu}

\author{Zsombor Moln\'{a}r}
\address{Institute of Applied Mechanics, University of Miskolc, Miskolc-Egyetemv\'{a}ros, H-3515, Hungary}
\ead{molnar.zsombor@student.uni-miskolc.hu}

\author{R\'{o}bert Kov\'{a}cs}
\address{Department of Energy Engineering, Faculty of Mechanical Engineering, Budapest University of Technology and Economics, M\H{u}egyetem rkp. 3., H-1111, Hungary \\ Department of Theoretical Physics, Wigner Research Centre for Physics, Budapest, Hungary}
\ead{kovacsrobert@energia.bme.hu}





\begin{abstract}
Thanks to the modern manufacturing technologies, heterogeneous materials with complex inner structure (e.g., foams) can be easily produced. However, their utilization is not straightforward as the classical constitutive laws are not necessarily valid. According to various experimental observations, the Guyer--Krumhansl equation stands as a promising candidate to model such complex structures. However, the practical applications need a reliable and efficient algorithm that is capable to handle both complex geometries and advanced heat equations. In the present paper, we present the development of a $hp$-type finite element technique, which can reliably applied. We investigate its convergence properties for various situations, being challenging in relation of stability and the treatment of fast propagation speeds. 
That algorithm is also proved to be outstandingly efficient, providing solutions four magnitudes faster than the commercial algorithms.

\end{abstract}


\end{frontmatter}


\section{Introduction}
The engineering practice utilizes several models to describe the material behaviour, these equations are called constitutive equations. The most frequent are the Navier--Stokes, Fourier and Hooke laws for continuum objects. The common point among them is that all define an equality:
\begin{align}
    q=-\lambda \partial_x T, \quad \Pi = - \nu \partial_x v, \quad \sigma = E \varepsilon, \label{eq1}
\end{align}
in which $q, \lambda$ and $T$ are the heat flux, thermal conductivity and temperature field; $\Pi$ describes the dynamic pressure as a consequences of the presence of velocity gradient $\partial_x v$; and the stress tensor $\sigma$ is proportional with the deformation $\varepsilon$. 
This simple structure usually allows to easily eliminate the current densities and prescribing the boundary conditions (BCs) using only one field variable, for instance, the temperature for heat conduction problems. It has also eased the development and implementation of various numerical and analytical solution techniques.

However, for advanced transport models, the situation becomes more difficult. The constitutive equation is not an equality but becomes a partial differential equation. Therefore, the conventional initial and BCs do not work without any further conditions \cite{Kov18gk, Kov22a}, and a more careful treatment is necessary. A glaring example is related to the so-called Guyer--Krumhansl (GK) equation,
\begin{equation}
    \tau\,\dot{q}+q+\lambda\,T'-\kappa^2\,q''=0\quad\mathrm{in}\;\Omega=[x|x\in(0,\ell)]\label{GKeq}
\end{equation}
in which the partial time derivative of the heat flux ($\dot q$) appears with a coefficient $\tau$ called relaxation time. Moreover, Eq.~\eqref{GKeq} consists of the second-order spatial derivative of $q$ as well, denoted by $q''$. Although the GK equation is usually known from the kinetic theory describing wave-type heat conduction phenomena in crystals \cite{GK66}, it also can be derived on a continuum basis \cite{VanFul12}. Thus the coefficient $\kappa^2$ is not related to the mean free path of phonons as in that continuum approach, the resulting model is independent of the micro-scale heat transfer mechanisms \cite{VanFul12}. Adding that the GK equation is successfully tested in heat conduction experiments on various rocks and foams as challenging problems \cite{Botetal16, FehKov21, LunEtal22}, this allows us to consider the GK equation as a promising candidate beyond Fourier and use it for more practical engineering applications. Hence it is essential to understand the properties of the GK equation and develop efficient and reliable numerical solution techniques. 

On that basis, a finite difference approach has been elaborated recently \cite{RietEtal18}, but it is strongly limited to simple (regular) geometries. That numerical methodology is also adapted for rheological and nonlinear models \cite{PozsEtal20}, with carefully investigating the role of initial and boundary conditions analytically, too \cite{Kov18gk, Kov22a}. 
Therefore, we aim to develop a $hp$-version finite element method (FEM), which is able to preserve the advantageous properties of the earlier finite difference scheme to handle the initial and boundary conditions reliably, and, also, to be adaptable for complex geometries and possess high accuracy and fast convergence properties.

The simplest classical FE techniques use linear piecewise polynomials to approximate the solution and achieve the desired accuracy with mesh refinement. The philosophy of using low-order polynomials over successively finer meshes is called $h$-type approximation technique \cite{Bat96,Ciar78,Hug87,Redd02}. In most cases \cite{DekaDutta19, YangEtal19, NazmEtal21, XuLi03}, that $h$-type approach is utilised for the Maxwell--Cattaneo--Vernotte- (MCV, with $\kappa^2=0$ in Eq.~\eqref{GKeq}), or for the dual-phase-lag (DPL) equations, which models have much less practical interest.

However, the conventional $h$-FEMs can provide slow convergences and low accuracy for specific types of initial and BCs, as well as when certain material and/or geometric parameters of the considered model problem is close to their limit value. In order to overcome these numerical difficulties, the $p$-version FE technology will be used as alternative, promising strategy, being originally-introduced in \cite{SzabBabu2011,BabSzabKatz81}. The idea is to keep the coarse mesh fixed, and the convergence and the high accuracy are achieved solely by increasing the polynomial degree $p$ of the approximated variables. It was proven that the rate of convergence is
much higher for $p$-FEM than that is possible with $h$-FEM, and exponential for smooth solutions and even for non-smooth solutions using properly chosen mesh refinement coupled with the increase of the polynomial degree $p$ ($hp$-strategy) \cite{BabSzabKatz81}. This is the main motivation behind our research, thus the $p$-type approximation technique is chosen to develop new FEM for the above-mentioned, refined theories of thermodynamics.

\subsection*{Problem formulation}
The GK equation \eqref{GKeq} is accompanied by the balance of internal energy for heat conduction processes, that is,
\begin{equation}
    \rho\,c_V\,\dot{T}+q'=0\quad\mathrm{in}\;\Omega=[x|x\in(0,\ell)], \label{bale}
\end{equation}
where $\rho$ and $c_V$ are the mass density and specific heat; the source terms are neglected. In the following, we present the discretization procedure for that system of equations.

For the sake of generality the system of basic differential Eqs. \eqref{GKeq}--\eqref{bale} is subjected the spatial descriptions 
\begin{align}
T(t,0)=\tilde{T}(t)\:,\label{Dirich}\\
q(t,\ell)=\tilde{q}(t)\label{Neu}
\end{align}as Dirichlet- and Neumann-type BCs at $x=0$ and $x=\ell$, respectively, as well as the temporal prescriptions
\begin{equation}
T(0,x)=T_0(x)\quad\mathrm{and}\quad q(0,x)=q_0(x)\label{init_cond}\end{equation}
as initial conditions at $t_0=0$ s.

\section{Two-field mixed $hp$-version finite element method}
In this section, new $hp$-FEMs will be presented for the numerical solution of the MCV- and the GK model which are based on two-field mixed variational formulations, see some similar procedures in \cite{Toth2018,Toth2016,TOTH2021}.    
\subsection{Variational formulations}
\subsubsection{Maxwell--Cattaneo--Vernotte model}
In this model problem, after testing Eq. \eqref{GKeq}, i.e., multiplying this with the test function $v$ and integrating it over the domain $\Omega$ with setting $\kappa^2=0$, then relaxing Eq. \eqref{bale},
i.e., multiplying this with the test function $u$ and
integrating it by parts, as well as building-in the BC \eqref{Neu} and
the homogeneous form of the BC \eqref{Dirich} for $u$, we arrived at a two-field mixed variational formulation which reads as follows.
Find the duet $T\in H^1(\Omega)$ and $q\in L^2(\Omega)$ as trial functions satisfying {\it a priori} the BC \eqref{Dirich} and the initial condition \eqref{init_cond}
such that
\begin{align}
-\int_{\Omega}\rho\,c_V\,\dot{T}\,u\,\mathrm{d}\Omega+\int_{\Omega}q\,u'\,\mathrm{d}\Omega&=\tilde{q}_{\ell}\,u(\ell)\quad\forall u\in H^1(\Omega)\:,\label{MCVvareq1}\\
\int_{\Omega}\tau\,\dot{q}\,v\,\mathrm{d}\Omega+\int_{\Omega}q\,v\,\mathrm{d}\Omega+\int_{\Omega}\lambda\,T'\,v\,\mathrm{d}\Omega&=0\quad\forall v\in L^2(\Omega)\label{MCVvareq2}
\end{align}satisfying {\it a priori} the homogeneous form of the BC \eqref{Dirich}. Here, $H^1(\Omega)$ stands for the Sobolev space of order 1 \cite{BrenScott08} and represents the regularity property for
$T$ and $u$ while $L^2(\Omega)$ defines square integrable function space for $q$ and $v$. 
\subsubsection{Guyer--Krumhansl model}
In this model problem relaxing again Eq. \eqref{bale} with the test function $u$, but now, after having tested Eq. \eqref{GKeq} with test function $v$, integrating-by-parts the last term of this variational equation, we obtain again a two-field mixed variational formulation but with a bit different regularity property from Eq. \eqref{MCVvareq2}. Accordingly, we seek the variable pair $(T,\,q)\in H^1(\Omega)$ satisfying {\it a priori} the BC \eqref{Dirich} and the initial condition \eqref{init_cond} in such a way that
\begin{align}
-\int_{\Omega}\rho\,c_V\,\dot{T}\,u\,\mathrm{d}\Omega+\int_{\Omega}q\,u'\,\mathrm{d}\Omega&=\tilde{q}_{\ell}\,u(\ell)\quad\forall u\in H^1(\Omega)\:,\label{GKvareq1}\\
\int_{\Omega}\tau\,\dot{q}\,v\,\mathrm{d}\Omega+\int_{\Omega}q\,v\,\mathrm{d}\Omega+\int_{\Omega}\lambda\,T'\,v\mathrm{d}\Omega+\int_{\Omega}\kappa^2\,q'\,v'\,\mathrm{d}\Omega
&=q'(\ell)\,v(\ell)-q'(0)\,v(0)\quad\forall v\in H^1(\Omega)\label{GKvareq2}
\end{align} hold true, ensuring {\it a priori} the homogeneous form of the BC \eqref{Dirich} once again.  Notably, again, for both variables $(T,\,q)\in H^1(\Omega)$,  which stands as a crucial difference for other FE approaches for heat equations.

\subsection{$hp$-finite element discretizations}
Let us consider now the $hp$-FE discretization of the domain $\Omega$. Thenceforward $\Omega$ is divided into $n$ physical sub-domain.
Then the master element $\Omega_{\mathrm{mas}}:=\{\,\eta\,|\,\eta\in(-1,1)\}$ is mapped onto the $i$-th physical element $K:=\{x\,|x^i\in(x_i,x_{i+1})\}\subset\Omega_{\mathrm{h}}$ with the nodal points $x_i$ and $x_{i+1}$ of the FE mesh $\Omega_{\mathrm{h}}:0<x_1<x_2<...<x_i<x_{i+1}<...<x_n<x_{n+1}=\ell$ by the mapping function
$x^i:=N_1(\eta)x_i+N_2(\eta)x_{i+1}$, where $N_1=(1-\eta)/2$ and $N_2=(1+\eta)/2$, $i=1,\ldots,n$.

The approximation spaces for the trial- and test function group, ($T,\,q$) and ($u,\,v$) are $hp$-FE function spaces consisting of piece-wise continuous polynomial functions, being spanned over the $i$-th element by the external shape functions $N_1$ and $N_2$ for $p=1$, as well as their supplemented set with the bubble modes $N_k(\eta)=[L_{k-1}(\eta)-L_{k-3}(\eta)]/\sqrt{2(2k-3)}$ for $p\geq2$, where $L_k(\eta)$ are the orthogonal Legendre polynomials, $k=3,4,\ldots,p+1$ \cite{SzabBabu2011,hpbook,DueRankSzab2017}. The independent trial- and test functions, $T,\;u$ and $q,\;v$ are approximated on the $i$-th element of $\Omega_{\mathrm{h}}$ by the polynomial degrees $p+1$ and $p$, or $p+1$, respectively, for the MCV- and the GK model according to Table \ref{tab:polynomial_spaces}, keeping in mind the regularity assumptions on the trial functions $T$ and $q$, as well as the test functions $u$ and $v$ in Eqs. \eqref{MCVvareq1}--\eqref{MCVvareq2} for the MCV model and in Eqs. \eqref{GKvareq1}--\eqref{GKvareq2} for the GK model. Here $p$ is the actual polynomial degree set to each element.
\begin{table}[ht]
\centering
\begin{tabular}{|c|c|c|} \hline
\backslashbox[2.5cm]{Model}{Function}
 & $T,\,u$ & $q,\,v$ \\\hline
MCV  & $p+1$ & $p$ \\\hline
GK  & $p+1$ & $p+1$ \\\hline
\end{tabular}\caption{Polynomial approximation spaces.}
\label{tab:polynomial_spaces}
\end{table}

Let us represent now the variational equations \eqref{MCVvareq1}--\eqref{MCVvareq2} and \eqref{GKvareq1}--\eqref{GKvareq2} in matrix form. 
For the sake of simplicity let $^{\mathrm{h}}\boldsymbol{t}$ and $^{\mathrm{h}}\boldsymbol{q}$ indicate
the column matrices corresponding to the time-dependent unknown coefficients of the temperature and the heat flow, respectively.
After carrying out the numerical integrations on each element $i$, on the spatial assembling process we have the following time-dependent matrix equation system
\begin{equation}\begin{array}{rl}
\mathcal{C}\:^{\mathrm{h}}\dot{\boldsymbol{t}}+\mathcal{Q}^T\:^{\mathrm{h}}\boldsymbol{q}&=\boldsymbol{d}\:,\\
\mathcal{T}\:^{\mathrm{h}}\dot{\boldsymbol{q}}+\mathcal{Q}\:^{\mathrm{h}}\boldsymbol{t}+\mathcal{K}\:^{\mathrm{h}}\boldsymbol{q}&=\boldsymbol{0}\:,
\end{array}\label{eq:time-dependent_matrix-equation}
\end{equation}
which can be written in the simplified form 
\begin{equation}
\mathcal{A}\:^{\mathrm{h}}\dot{\boldsymbol{\alpha}}+\mathcal{B}\:^{\mathrm{h}}\boldsymbol{\alpha}=\boldsymbol{f}\:,\label{system}    
\end{equation}
where
\begin{equation}^{\mathrm{h}}\boldsymbol{\alpha}=\left[\begin{array}{c}
    ^{\mathrm{h}}\boldsymbol{t}\\
    ^{\mathrm{h}}\boldsymbol{q}
    \end{array}\right]\:,\quad
    \mathcal{A}=\left[\begin{array}{cc}
    \mathcal{C}&\boldsymbol{0}\\
    \boldsymbol{0}&\mathcal{T}
    \end{array}\right]\;,\quad\mathcal{B}=\left[\begin{array}{cc}
    \boldsymbol{0}&\mathcal{Q}^T\\
    \mathcal{Q}&\mathcal{K}
    \end{array}\right]\quad\mathrm{and}\quad\boldsymbol{f}=\left[\begin{array}{c}
    \boldsymbol{d}\\
    \boldsymbol{0}
    \end{array}\right]\;,
\end{equation}
in which $\mathcal{C}$, $\mathcal{T}$ and $\mathcal{K}$ denote the consistent matrices of the specific heat, the relaxation term and the heat conductivity, respectively,
while $\mathcal{Q}$ is the consistent coupling matrix of the system.

\section{Numerical experiments}

In this section the newly-developed two-field $hp$-version FEMs will be tested on a representative initial-boundary value problem,
focusing on the transient analyzes. In order to demonstrate the computational performance of the $hp$-type mixed FEs, first we investigate the point-wise $h$- and $p$-convergence
of the relative error that are measured, respectively, in the maximum norm of the temperature and the heat flow
\begin{equation}
e^{\mathrm{hp}}=\frac{\underset{t\in[t_0,t_1]}{\max}\left|\:^{\mathrm{h}}\boldsymbol{t}-\boldsymbol{t}_{\mathrm{ref}}\right|}
{\underset{t\in[t_0,t_1]}{\max}\left|\boldsymbol{t}_{\mathrm{ref}}\right|}\quad\textrm{and}\quad
e^{\mathrm{hp}}=\frac{\underset{t\in[t_0,t_1]}{\max}\left|\:^{\mathrm{h}}\boldsymbol{q}-\boldsymbol{q}_{\mathrm{ref}}\right|}
{\underset{t\in[t_0,t_1]}{\max}\left|\boldsymbol{q}_{\mathrm{ref}}\right|}\;.\label{error_measures}
\end{equation}
The related thermodynamic system \eqref{system} is solved numerically for the relatively small time interval $t\in[0,10]$ s, using the implicit time integration scheme \cite{Hug87}.
The number of the constant time step is set to $n_t=10000$, yielding the time step size $\Delta t=0.001$ s. The considered bodies of length $\ell=0.005$ m are made of rock-like materials which have $c_V=800\;\mathrm{(J/kgK)}$, $\rho=2600\;\mathrm{kg/m^{3}}$. Additionally, the new parameters $\tau$ and $\kappa^2$ are determined based on earlier experimental data \cite{Botetal16} to be realistic, therefore we choose $\tau=\{0.05, 0.15, 0.3\}$ s, with $\kappa^2=8\cdot 10^{-6}$ m$^2$.

However, in order to challenge the numerical procedure, we also test for unrealistically high $\kappa^2=0.8$ m$^2$ as that coefficient greatly influences the characteristic speed of diffusion and makes the propagation faster with five magnitudes with respect to the previous situations.

Furthermore, these are subjected to the time-dependent Neumann-type BCs
\begin{equation}\tilde{q}_0(t)=10000\,\frac{c_1\,c_2}{c_2-c_1}
\left[\exp\left(-c_1\frac{ t}{t_p}\right)-\exp\left(-c_2\frac{ t}{t_p}\right)\right]\quad\mathrm{and}\quad\tilde{q}_{\ell}(t)=0\end{equation}
as prescribed heat flows and the initial conditions $T_0(x)=293$ K and $q_0(x)=0$, where $c_1=1/0.075$, $c_2=6$, $t_p=0.008$ s, namely, these are initially at rest, i.e, in equilibrium state. The $c_1$ and $c_2$ coefficients are chosen based on \cite{FehKov21}, to keep the excitation to be realistic as much as possible.

When the computational performance of the $p$-approximation is analyzed, a 8-, 52- and 20-element uniform mesh kept fixed on the domain $\Omega$
for the GK model with the parameter settings $\kappa^2=0.8$ m$^2$, $\kappa^2=0.000008$ m$^2$ and the MCV model, respectively,
while the polynomial degree $p$ is ranging from 2 to 8 for each element. During the $h$-convergence studies,
the domain $\Omega$ is uniformly-refined in 7 steps from the element number $n=8$, 52 and 20 to 20, 88 and 44, respectively,
for the GK model with $\kappa^2=0.8$ m$^2$, $\kappa^2=0.000008$ m$^2$ and the MCV model. In each step the polynomial degree $p$ is
set to 2 and kept fixed on each element. At each step, the number of degrees of freedom (DOF) is calculated as the sum of
the total unknown coefficients occurring in $\boldsymbol{\alpha}$.

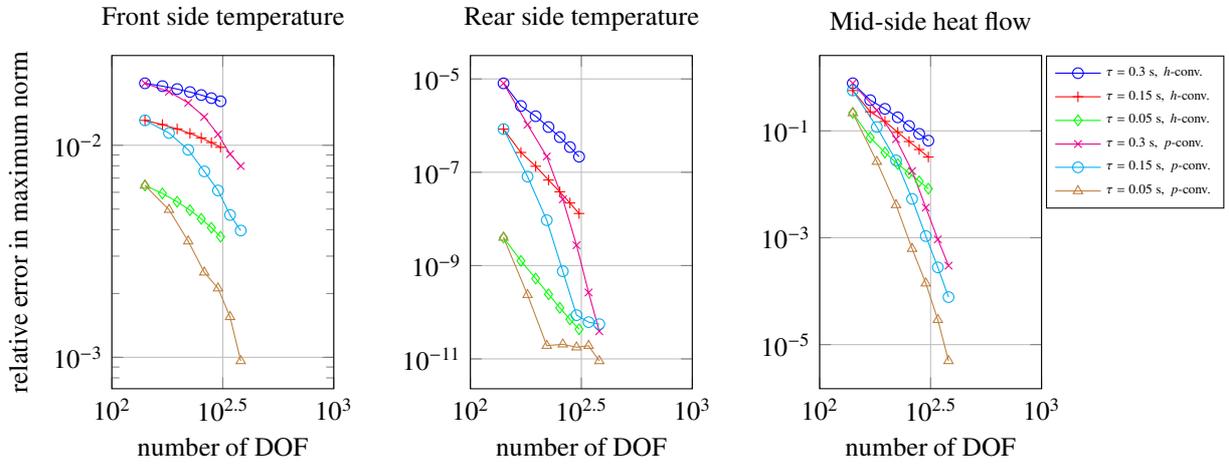
\begin{figure}[H]
\begin{tabular}{lll}
\begin{tikzpicture}
\tikzset{every mark/.append style={scale=1}}
\begin{loglogaxis}[grid=major,width=4.5cm,height=6cm,xmin=100,xmax=1000,
legend style={font=\tiny},legend cell align=left,
xlabel=number of DOF,ylabel=relative error in maximum norm,title={Front side temperature}]
\addplot[color=blue,mark=o] table[x=dh,y=h03] {mcvkonvfrontside.dat};
\addplot[color=red,mark=+] table[x=dh,y=h015] {mcvkonvfrontside.dat};
\addplot[color=green,mark=diamond] table[x=dh,y=h005] {mcvkonvfrontside.dat};
\addplot[color=magenta,mark=x] table[x=dp,y=p03] {mcvkonvfrontside.dat};
\addplot[color=cyan,mark=o] table[x=dp,y=p015] {mcvkonvfrontside.dat};
\addplot[color=brown,mark=triangle] table[x=dp,y=p005] {mcvkonvfrontside.dat};
\end{loglogaxis}
\end{tikzpicture}&
\begin{tikzpicture}
\tikzset{every mark/.append style={scale=1}}
\begin{loglogaxis}[grid=major,width=4.5cm,height=6cm,xmin=100,xmax=1000,
legend style={font=\tiny},legend cell align=left,
xlabel=number of DOF,title={Rear side temperature}]
\addplot[color=blue,mark=o] table[x=dh,y=h03] {mcvkonvrearside.dat};
\addplot[color=red,mark=+] table[x=dh,y=h015] {mcvkonvrearside.dat};
\addplot[color=green,mark=diamond] table[x=dh,y=h005] {mcvkonvrearside.dat};
\addplot[color=magenta,mark=x] table[x=dp,y=p03] {mcvkonvrearside.dat};
\addplot[color=cyan,mark=o] table[x=dp,y=p015] {mcvkonvrearside.dat};
\addplot[color=brown,mark=triangle] table[x=dp,y=p005] {mcvkonvrearside.dat};
\end{loglogaxis}\end{tikzpicture}&
\begin{tikzpicture}
\pgfplotsset{every axis legend/.append style={
at={(1.02,1)},
anchor=north west}}
\begin{loglogaxis}[grid=major,width=4.5cm,height=6cm,xmin=100,xmax=1000,
legend style={font=\tiny},legend cell align=left,
xlabel=number of DOF,title={Mid-side heat flow}]
\addplot[color=blue,mark=o] table[x=dh,y=h03] {mcvkonvmiddle.dat};
\addplot[color=red,mark=+] table[x=dh,y=h015] {mcvkonvmiddle.dat};
\addplot[color=green,mark=diamond] table[x=dh,y=h005] {mcvkonvmiddle.dat};
\addplot[color=magenta,mark=x] table[x=dp,y=p03] {mcvkonvmiddle.dat};
\addplot[color=cyan,mark=o] table[x=dp,y=p015] {mcvkonvmiddle.dat};
\addplot[color=brown,mark=triangle] table[x=dp,y=p005] {mcvkonvmiddle.dat};
\legend{{$\tau=0.3\;\mathrm{s},\;h$-conv.},{$\tau=0.15\;\mathrm{s},\;h$-conv.},{$\tau=0.05\;\mathrm{s},\;h$-conv.},
{$\tau=0.3\;\mathrm{s},\;p$-conv.},{{$\tau=0.15\;\mathrm{s},\;p$-conv.}},{$\tau=0.05\;\mathrm{s},\;p$-conv.}}
\end{loglogaxis}
\end{tikzpicture}
\end{tabular}
\caption{Convergence histories of the relative errors measured in maximum norm for the front- and rear side temperature, as well as the mid-side heat flow -- MCV model}
\label{fig:MCV}
\end{figure}

\begin{figure}[H]
\begin{tabular}{lll}
\begin{tikzpicture}
\tikzset{every mark/.append style={scale=1}}
\begin{loglogaxis}[grid=major,width=4.5cm,height=6cm,xmin=200,xmax=2000,
xlabel=number of DOF,ylabel=relative error in maximum norm,title={Front side temperature}]
\addplot[color=blue,mark=o] table[x=dh,y=h03] {gk0000008konvfrontside.dat};
\addplot[color=red,mark=+] table[x=dh,y=h015] {gk0000008konvfrontside.dat};
\addplot[color=green,mark=diamond] table[x=dh,y=h005] {gk0000008konvfrontside.dat};
\addplot[color=magenta,mark=x] table[x=dp,y=p03] {gk0000008konvfrontside.dat};
\addplot[color=cyan,mark=o] table[x=dp,y=p015] {gk0000008konvfrontside.dat};
\addplot[color=brown,mark=triangle] table[x=dp,y=p005] {gk0000008konvfrontside.dat};
\end{loglogaxis}
\end{tikzpicture}&
\begin{tikzpicture}
\tikzset{every mark/.append style={scale=1}}
\begin{loglogaxis}[grid=major,width=4.5cm,height=6cm,xmin=200,xmax=2000,
legend style={font=\tiny},legend cell align=left,
xlabel=number of DOF,title={Rear side temperature}]
\addplot[color=blue,mark=o] table[x=dh,y=h03] {gk0000008konvrearside.dat};
\addplot[color=red,mark=+] table[x=dh,y=h015] {gk0000008konvrearside.dat};
\addplot[color=green,mark=diamond] table[x=dh,y=h005] {gk0000008konvrearside.dat};
\addplot[color=magenta,mark=x] table[x=dp,y=p03] {gk0000008konvrearside.dat};
\addplot[color=cyan,mark=o] table[x=dp,y=p015] {gk0000008konvrearside.dat};
\addplot[color=brown,mark=triangle] table[x=dp,y=p005] {gk0000008konvrearside.dat};
\end{loglogaxis}\end{tikzpicture}&
\begin{tikzpicture}
\pgfplotsset{every axis legend/.append style={
at={(1.02,1)},
anchor=north west}}
\begin{loglogaxis}[grid=major,width=4.5cm,height=6cm,xmin=200,xmax=2000,
legend style={font=\tiny},legend cell align=left,
xlabel=number of DOF,title={Mid-side heat flow}]
\addplot[color=blue,mark=o] table[x=dh,y=h03] {gk0000008konvmiddle.dat};
\addplot[color=red,mark=+] table[x=dh,y=h015] {gk0000008konvmiddle.dat};
\addplot[color=green,mark=diamond] table[x=dh,y=h005] {gk0000008konvmiddle.dat};
\addplot[color=magenta,mark=x] table[x=dp,y=p03] {gk0000008konvmiddle.dat};
\addplot[color=cyan,mark=o] table[x=dp,y=p015] {gk0000008konvmiddle.dat};
\addplot[color=brown,mark=triangle] table[x=dp,y=p005] {gk0000008konvmiddle.dat};
\legend{{$\tau=0.3\;\mathrm{s},\;h$-conv.},{$\tau=0.15\;\mathrm{s},\;h$-conv.},{$\tau=0.05\;\mathrm{s},\;h$-conv.},
{$\tau=0.3\;\mathrm{s},\;p$-conv.},{{$\tau=0.15\;\mathrm{s},\;p$-conv.}},{$\tau=0.05\;\mathrm{s},\;p$-conv.}}
\end{loglogaxis}
\end{tikzpicture}
\end{tabular}
\caption{Convergence histories of the relative errors measured in maximum norm for the front- and rear side temperature, as well as the mid-side heat flow -- GK model with $\kappa^2=0.000008\;\mathrm{m}^2$}\label{fig:GK-small}
\end{figure}
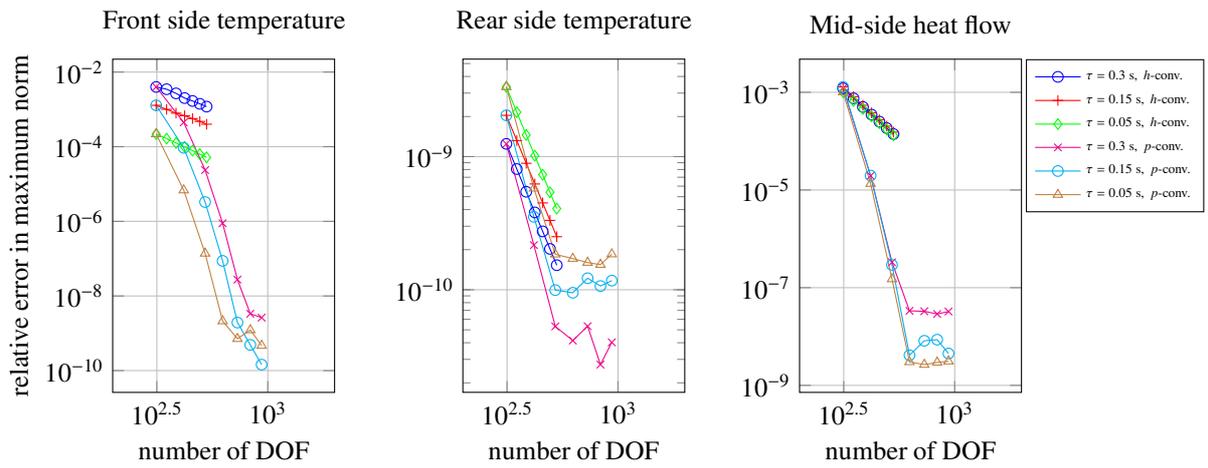

The relative error-convergence curves of the front- and rear side temperature, $T(t,0)$ and $T(t,\ell)$, as well as the mid-side heat flow $q(t,\ell/2)$ obtained for the MCV model and the GK model with relatively large and small value of $\kappa^2$ (0.8 m$^2$ and 0.000008 m$^2$) are plotted against the number of DOF on log-log scales in Figs. \ref{fig:MCV} and \ref{fig:GK-large}--\ref{fig:GK-small} for the relaxation time $\tau=0.3$ s, $\tau=0.15$ s and $\tau=0.05$ s, respectively.

Following from the characteristics of the relative error-curves, monotonically high rates are experienced both for
$h$- and for $p$-convergence. As expected, the $p$-convergence is much faster than the $h$-convergence. Namely,
the higher-order polynomial approximation helps a lot to increase the convergence rate. Thus, the desired accuracy is achieved with the use of much less number of DOF by the $p$-approximation. In the asymptotic range, the exponential type of the $p$-convergence behavior is observed while the $h$-convergence shows rather an algebraic type of convergence.

It can also be experienced that there is no significant influence of the relaxation time $\tau$ on the convergence rates. 
However, the convergence curves are shifted down (without changing their slopes) as $\tau$ is decreased achieving a lower relative error value. 
Besides, it can be seen that the computational limit is reached at a very small error level (at approximately $10^{-8}$).

\begin{figure}[H]
\centering\includegraphics[scale=0.6]{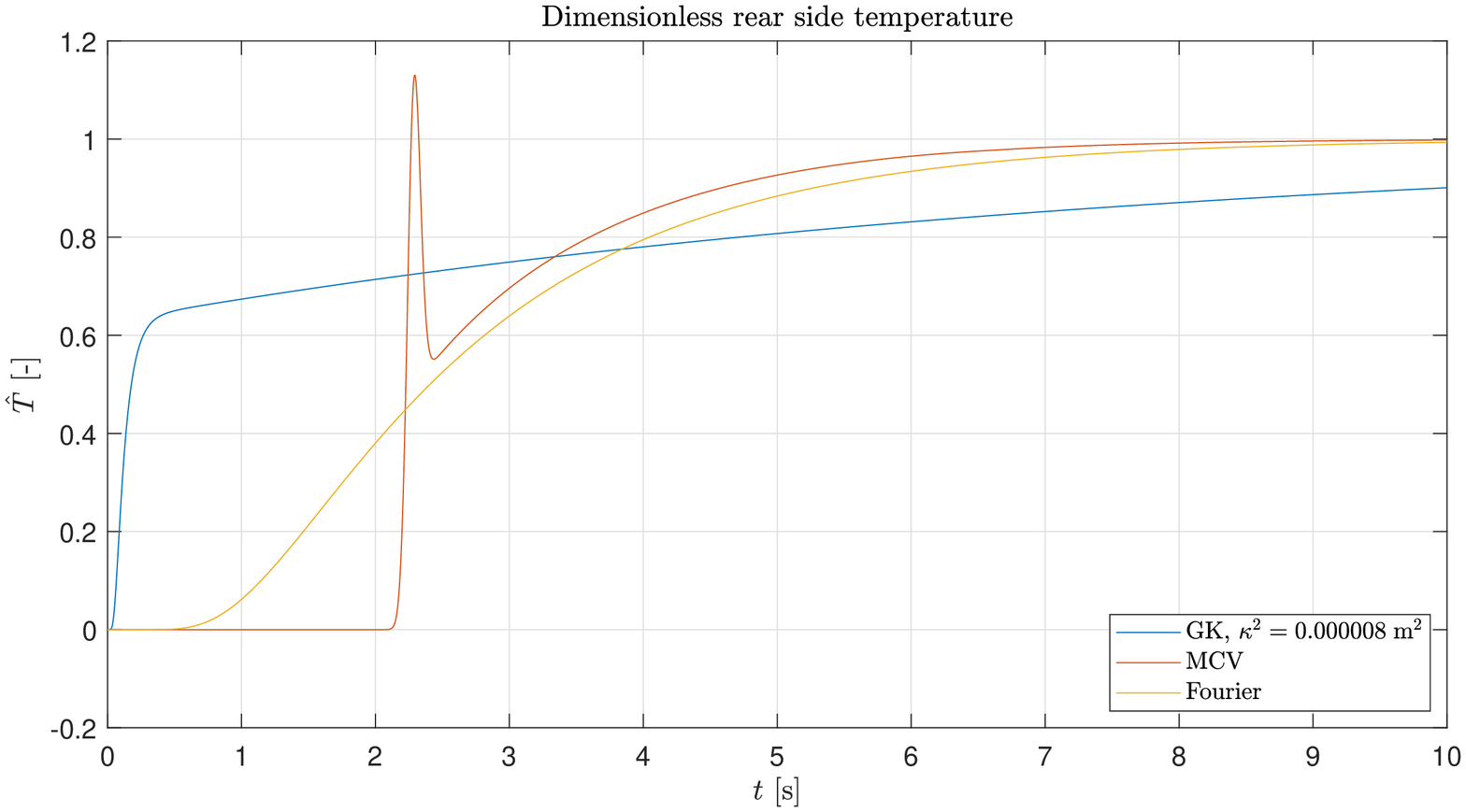}
\caption{Time series of the dimensionless rear side temperature for the MCV- and GK model with $\kappa^2=0.000008\;\mathrm{m}^2$ -- $\tau=0.3$ s, $p=10$ and $n=100$}\label{fig:T_rearside0.3}
\end{figure}

\begin{figure}[H]
\centering\includegraphics[scale=0.6]{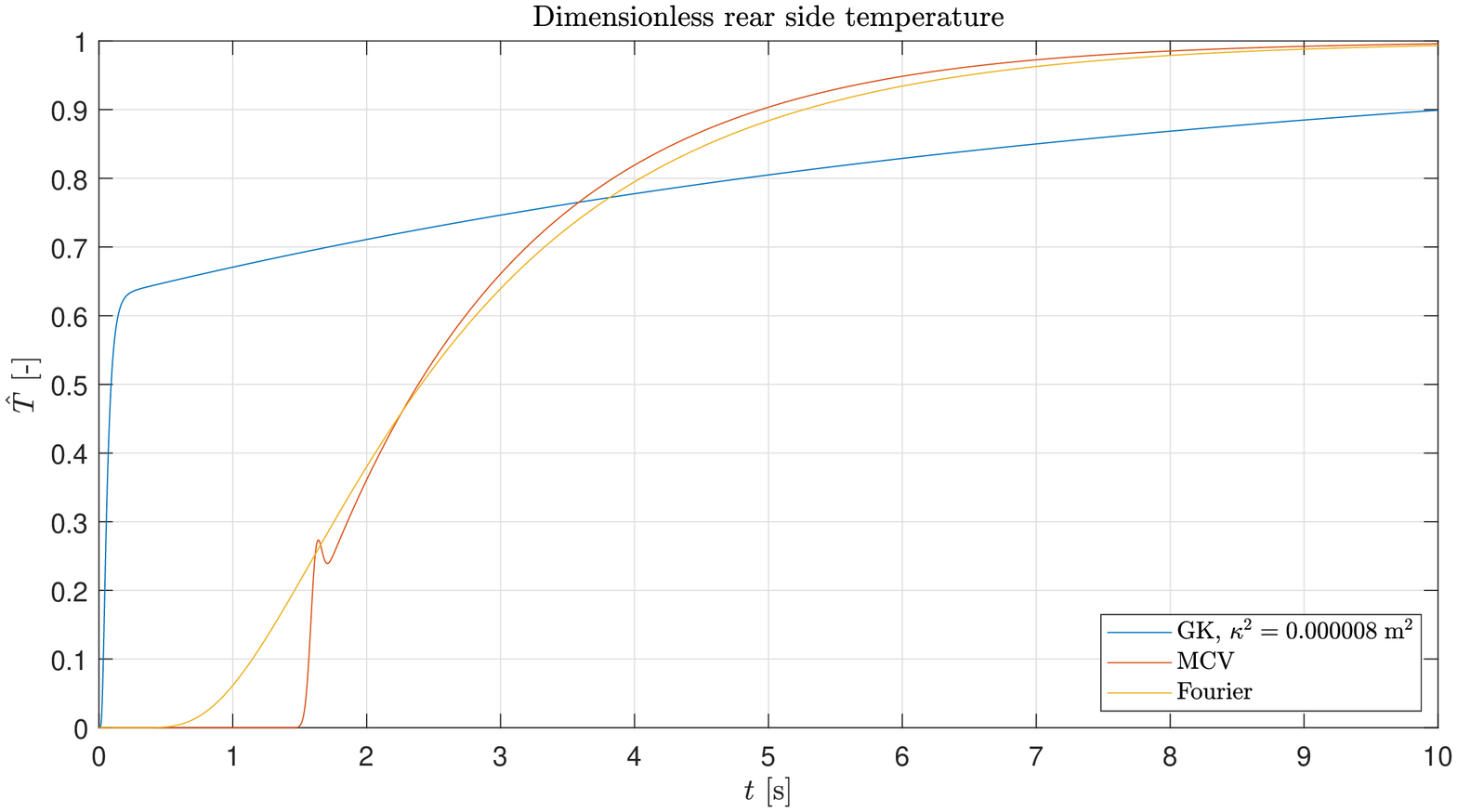}
\caption{Time series of the dimensionless rear side temperature for the MCV- and GK model with $\kappa^2=0.000008\;\mathrm{m}^2$ -- $\tau=0.15$ s, $p=10$ and $n=100$}\label{fig:T_rearside0.15}
\end{figure}

\begin{figure}[H]
\centering\includegraphics[scale=0.6]{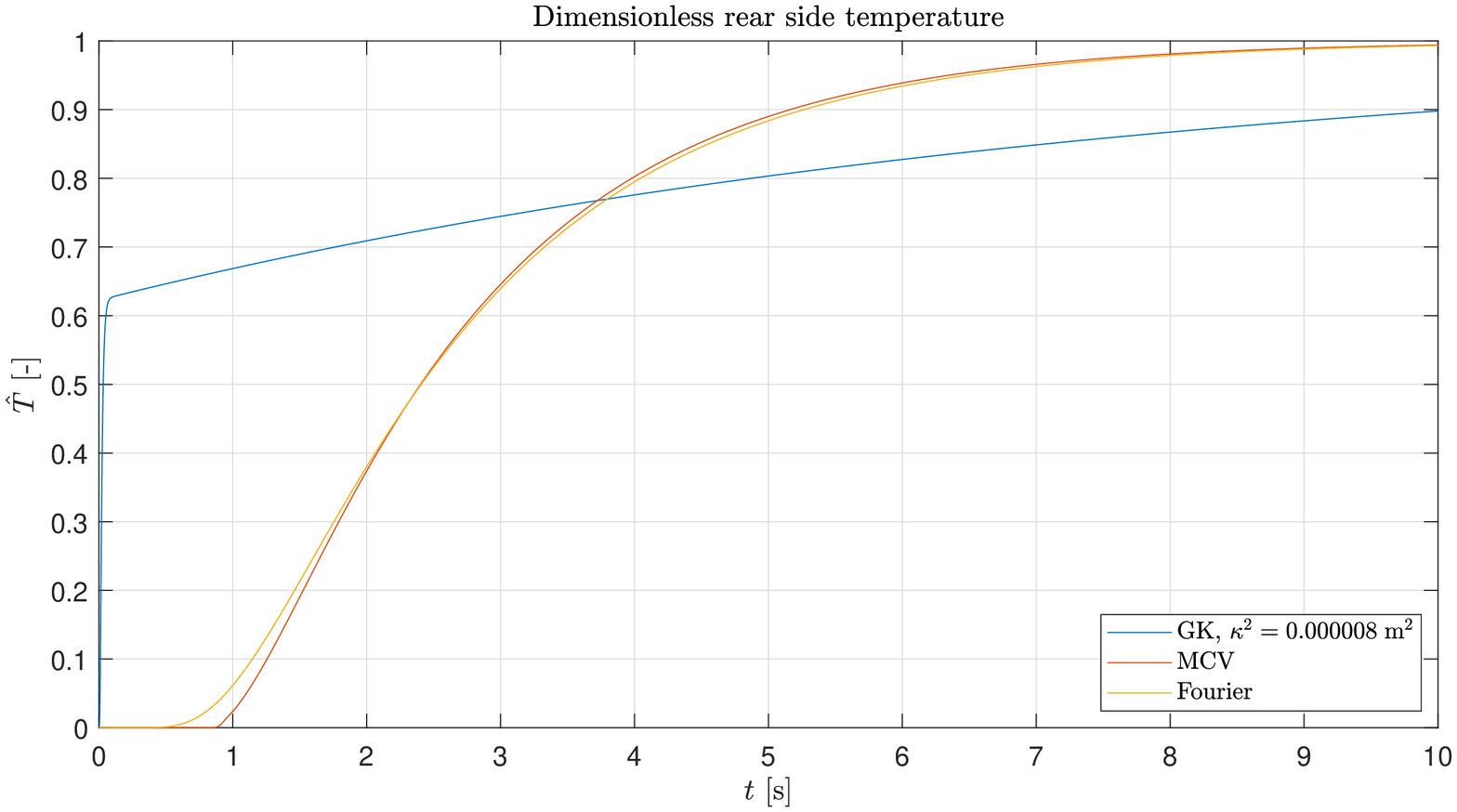}
\caption{Time series of the dimensionless rear side temperature for the MCV- and GK model with $\kappa^2=0.000008\;\mathrm{m}^2$ -- $\tau=0.05$ s, $p=10$ and $n=100$}\label{fig:T_rearside0.05}
\end{figure}

\begin{figure}[H]
\centering\includegraphics[scale=0.6]{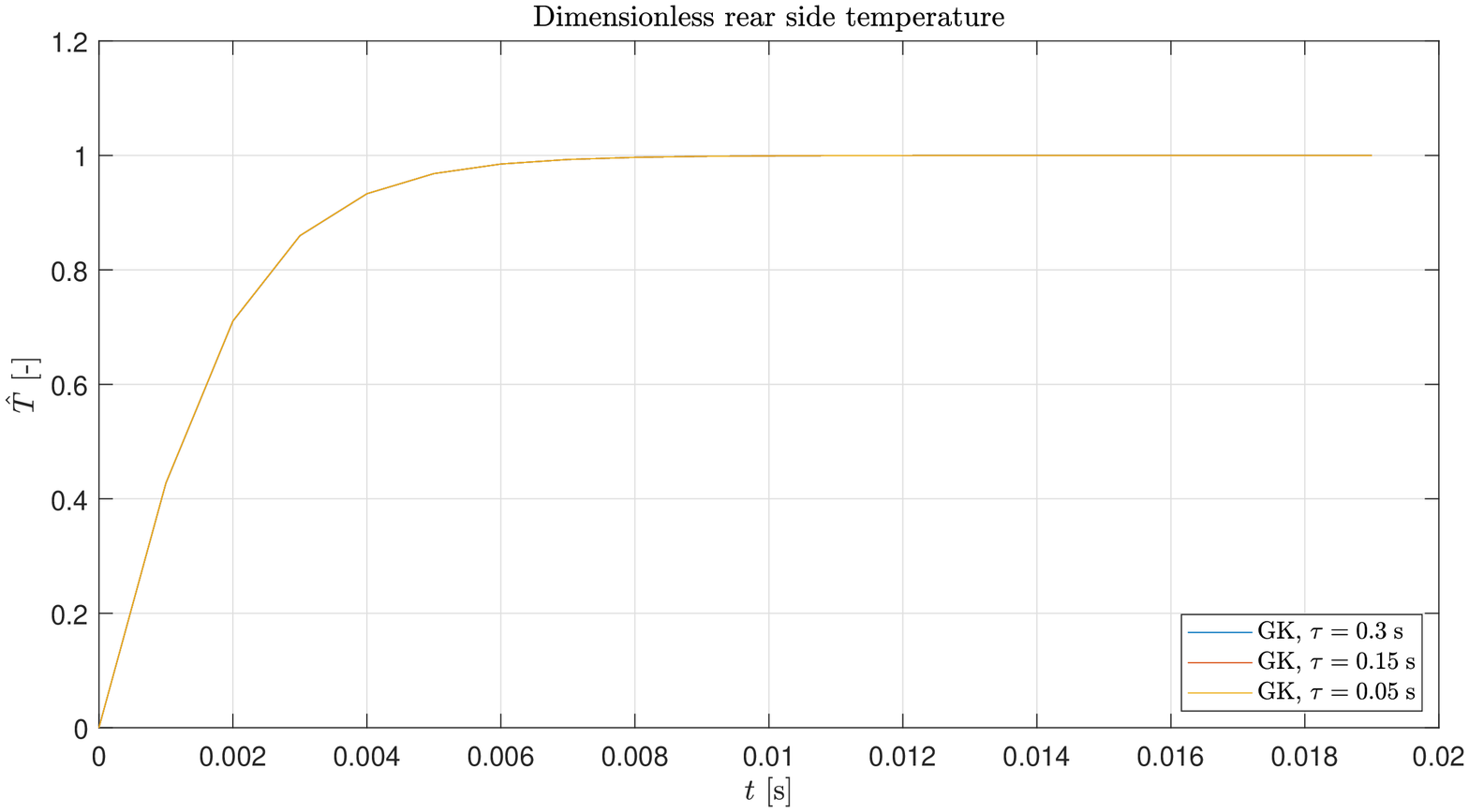}
\caption{Influence of $\tau$ on the time series of the dimensionless rear side temperature for the GK model with the over-diffuse setting $\kappa^2=0.8\;\mathrm{m}^2$ -- $p=10$ and $n=100$}\label{fig:T_rearside_kappa0.8}
\end{figure}

\begin{figure}[H]
\centering\includegraphics[scale=0.6]{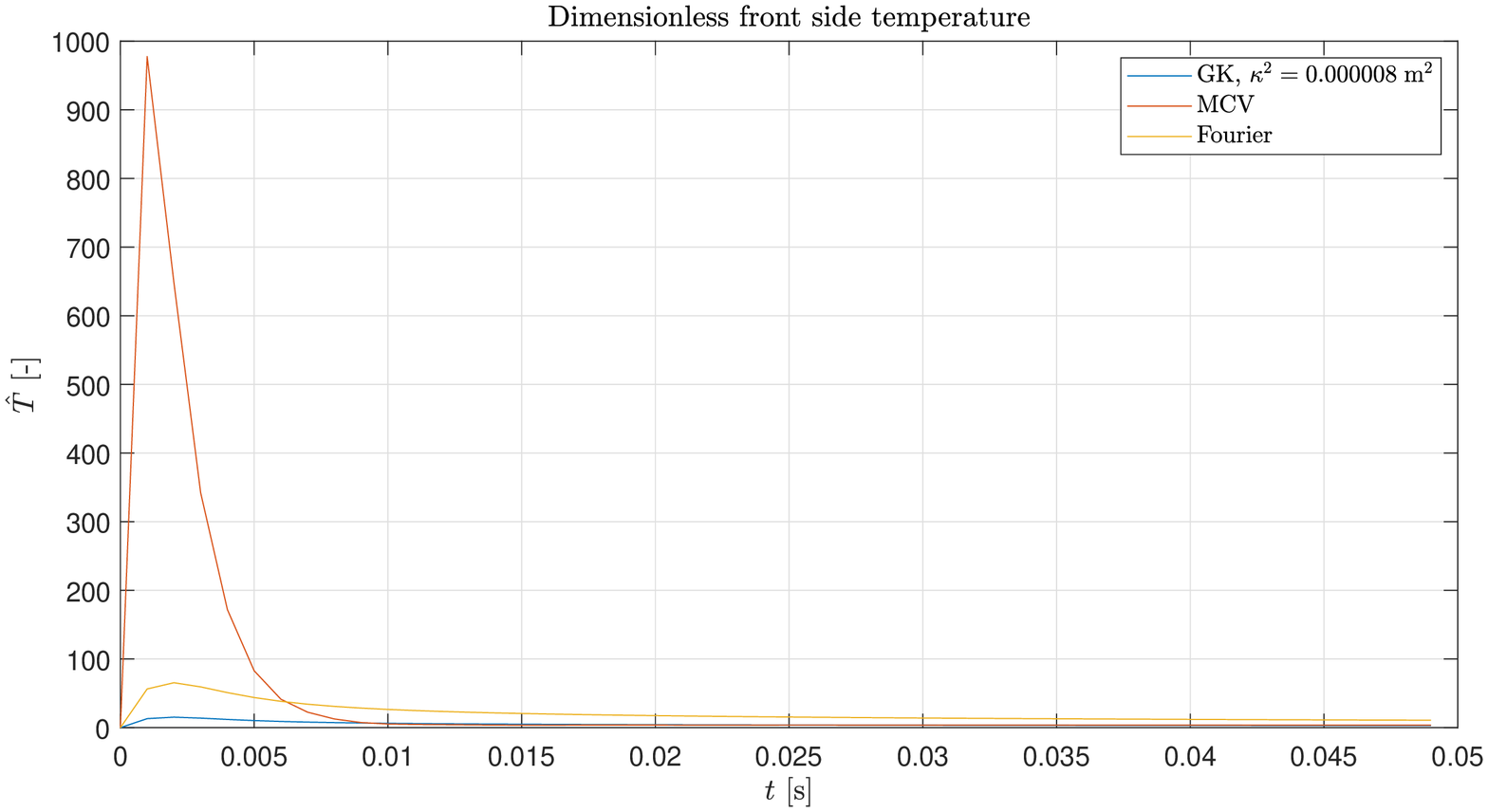}
\caption{Time series of the dimensionless front side temperature for the MCV- and GK model with $\kappa^2=0.000008\;\mathrm{m}^2$ -- $\tau=0.3$ s, $p=10$ and $n=100$}\label{fig:T_frontside0.3}
\end{figure}

\begin{figure}[H]
\centering\includegraphics[scale=0.6]{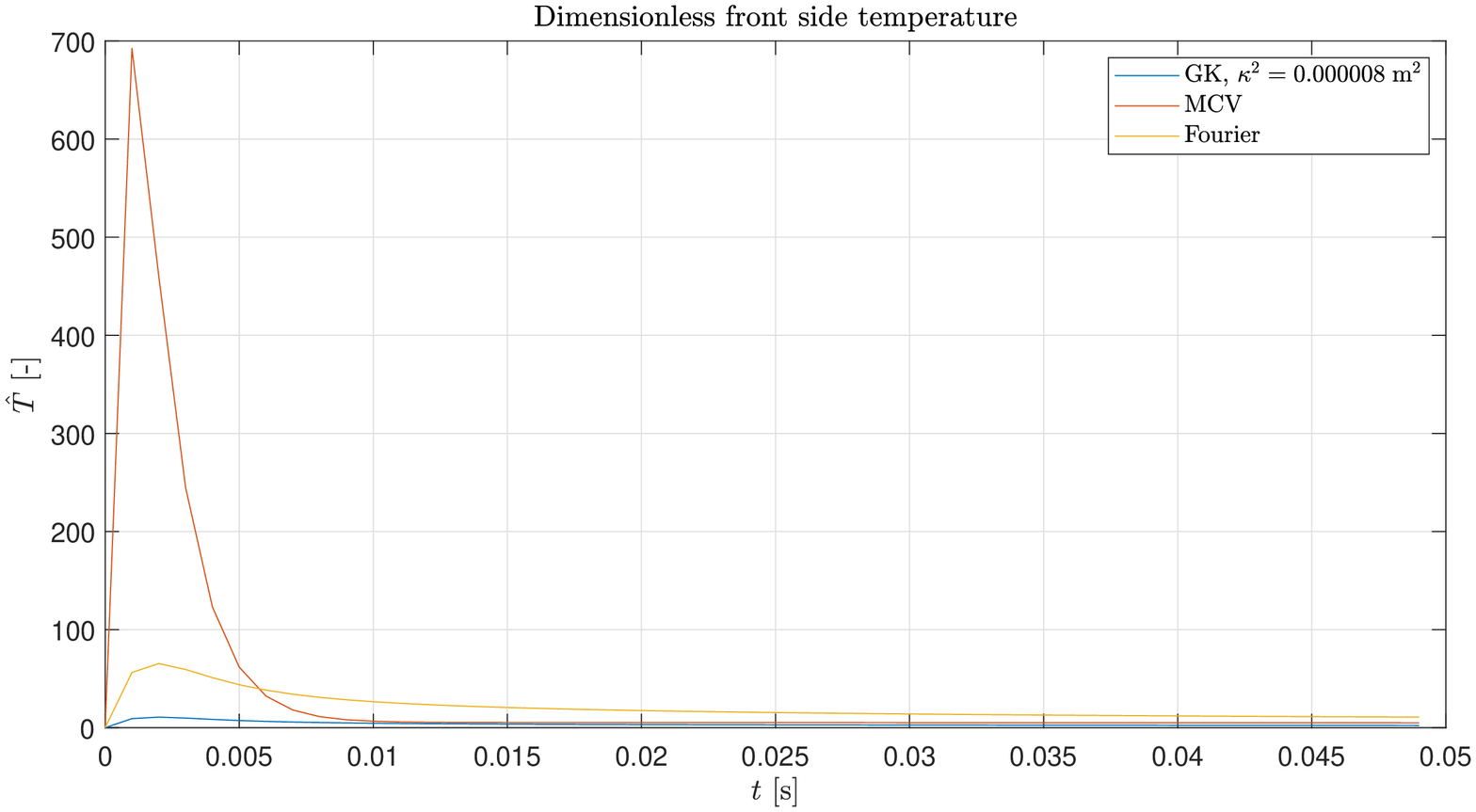}
\caption{Time series of the dimensionless front side temperature for the MCV- and GK model with $\kappa^2=0.000008\;\mathrm{m}^2$ -- $\tau=0.15$ s, $p=10$ and $n=100$}\label{fig:T_frontside0.15}
\end{figure}

\begin{figure}[H]
\centering\includegraphics[scale=0.6]{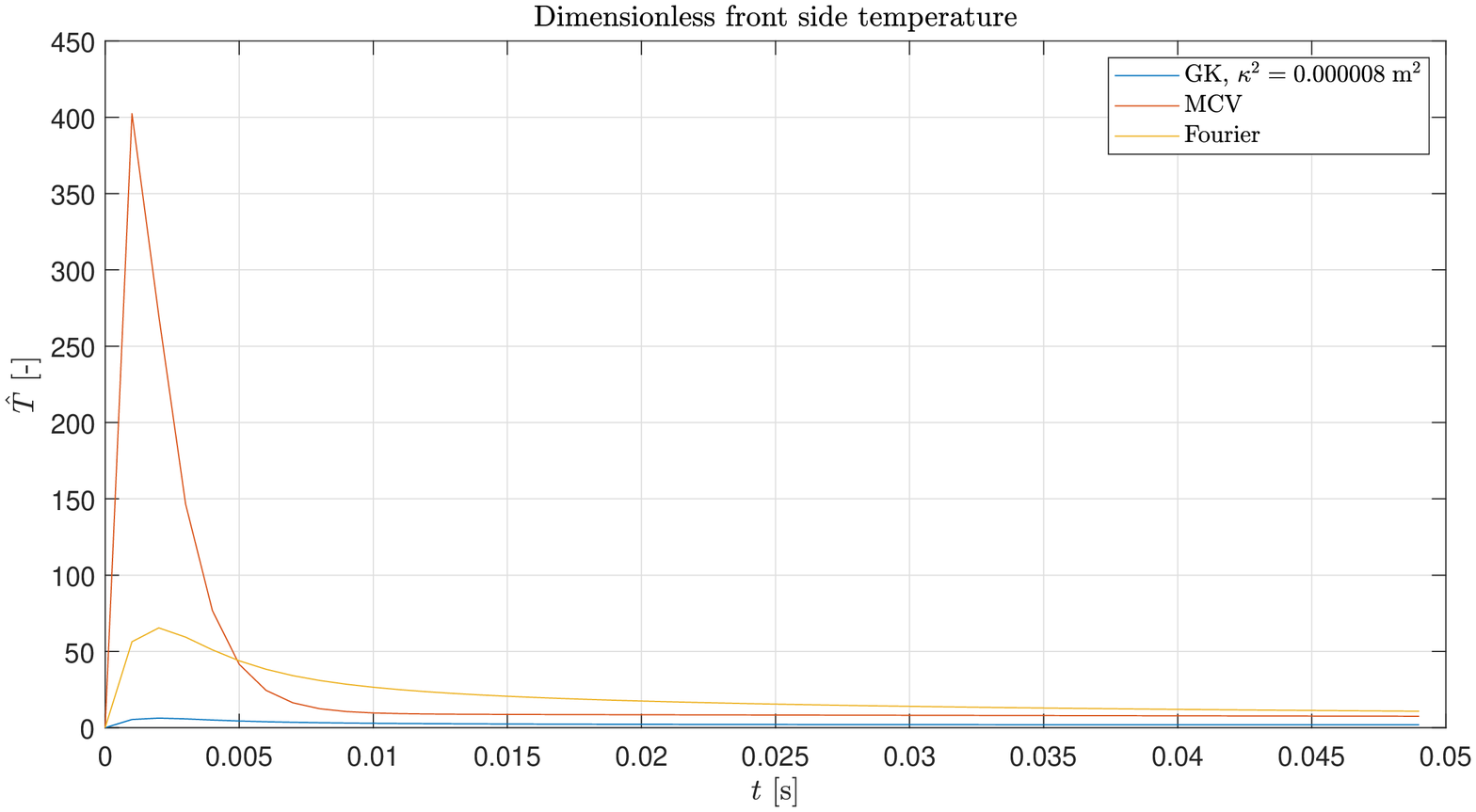}
\caption{Time series of the dimensionless front side temperature for the MCV- and GK model with $\kappa^2=0.000008\;\mathrm{m}^2$ -- $\tau=0.05$ s, $p=10$ and $n=100$}\label{fig:T_frontside0.05}
\end{figure}

\begin{figure}[H]
\centering\includegraphics[scale=0.6]{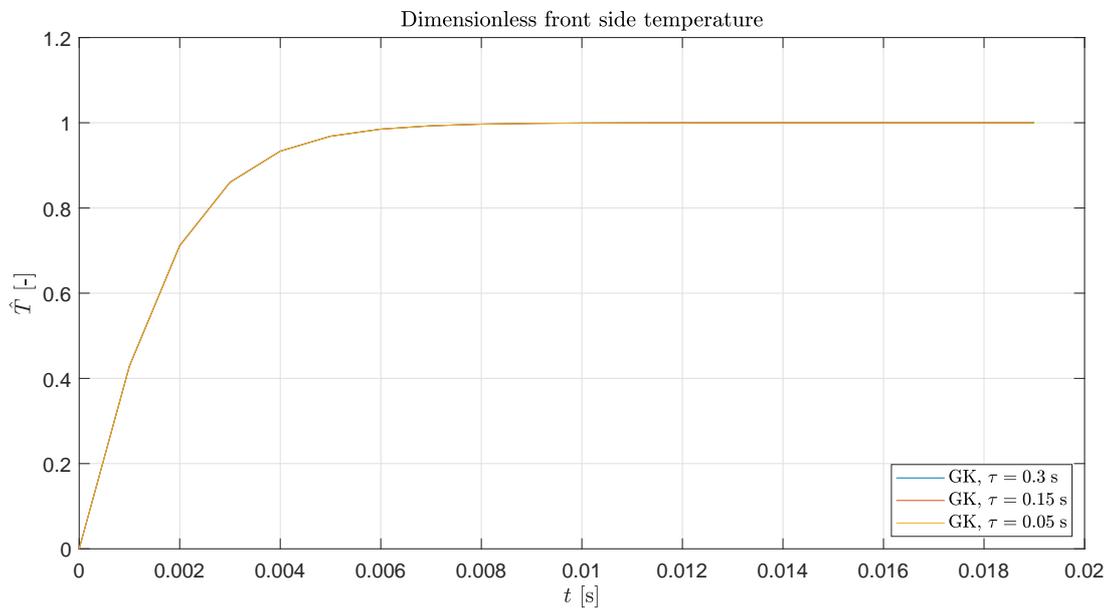}
\caption{Influence of $\tau$ on the time series of the dimensionless front side temperature for the GK model with the over-diffuse setting $\kappa^2=0.8\;\mathrm{m}^2$ -- $p=10$ and $n=100$}\label{fig:T_frontside_kappa08}
\end{figure}

Giving now as the illustration of the $hp$-FE solutions for the MCV model and the GK model with the quite small value of $\kappa^2$ (0.000008 m$^2$), the histories of the dimensionless rear- and front side temperatures are depicted, separately for three different relaxation times $\tau=0.3$ s, $\tau=0.15$ s and $\tau=0.05$ s. There is no any computational reason behind the dimensionless rescaling of temperature, this is used only for better overviewing the plots as the model is linear, and solely the characteristics are important. The adiabatic steady-state temperature is set to $1$ for these Figures. In Figs. \ref{fig:T_rearside0.3}--\ref{fig:T_rearside0.05} and, \ref{fig:T_frontside0.3}--\ref{fig:T_frontside0.05}, respectively, using a 100-element-mesh with the relatively high polynomial degree $p=10$ and comparing all the results to the $hp$-FE solution of the classical, Fourier model.  Figs. \ref{fig:T_rearside_kappa0.8} and \ref{fig:T_frontside_kappa08} represent the time series of the dimensionless rear- and front side temperature for \say{over-diffuse}  
thermodynamic system, i.e., for $\kappa^2=0.8$ m$^2$, zooming the $hp$-FE solutions in the rapidly change and very small time region $[0,0.02]$ s. These figures exhibit very well that the $\tau$-value has no significant effect on the solution for \say{over-diffuse} system.     

\begin{figure}[H]
\begin{tabular}{lll}
\begin{tikzpicture}
\tikzset{every mark/.append style={scale=1}}
\begin{loglogaxis}[grid=major,width=4.5cm,height=6cm,xmin=20,xmax=200,
xlabel=number of DOF,ylabel=relative error in maximum norm,title={Front side temperature}]
\addplot[color=blue,mark=o] table[x=dh,y=h03] {gk08konvfrontside.dat};
\addplot[color=red,mark=+] table[x=dh,y=h015] {gk08konvfrontside.dat};
\addplot[color=green,mark=diamond] table[x=dh,y=h005] {gk08konvfrontside.dat};
\addplot[color=magenta,mark=x] table[x=dp,y=p03] {gk08konvfrontside.dat};
\addplot[color=cyan,mark=o] table[x=dp,y=p015] {gk08konvfrontside.dat};
\addplot[color=brown,mark=triangle] table[x=dp,y=p005] {gk08konvfrontside.dat};
\end{loglogaxis}
\end{tikzpicture}&
\begin{tikzpicture}
\tikzset{every mark/.append style={scale=1}}
\begin{loglogaxis}[grid=major,width=4.5cm,height=6cm,xmin=20,xmax=200,
legend style={font=\tiny},legend cell align=left,
xlabel=number of DOF,title={Rear side temperature}]
\addplot[color=blue,mark=o] table[x=dh,y=h03] {gk08konvrearside.dat};
\addplot[color=red,mark=+] table[x=dh,y=h015] {gk08konvrearside.dat};
\addplot[color=green,mark=diamond] table[x=dh,y=h005] {gk08konvrearside.dat};
\addplot[color=magenta,mark=x] table[x=dp,y=p03] {gk08konvrearside.dat};
\addplot[color=cyan,mark=o] table[x=dp,y=p015] {gk08konvrearside.dat};
\addplot[color=brown,mark=triangle] table[x=dp,y=p005] {gk08konvrearside.dat};
\end{loglogaxis}\end{tikzpicture}&
\begin{tikzpicture}
\pgfplotsset{every axis legend/.append style={
at={(1.02,1)},
anchor=north west}}
\begin{loglogaxis}[grid=major,width=4.5cm,height=6cm,xmin=20,xmax=200,
legend style={font=\tiny},legend cell align=left,
xlabel=number of DOF,title={Mid-side heat flow}]
\addplot[color=blue,mark=o] table[x=dh,y=h03] {gk08konvmiddle.dat};
\addplot[color=red,mark=+] table[x=dh,y=h015] {gk08konvmiddle.dat};
\addplot[color=green,mark=diamond] table[x=dh,y=h005] {gk08konvmiddle.dat};
\addplot[color=magenta,mark=x] table[x=dp,y=p03] {gk08konvmiddle.dat};
\addplot[color=cyan,mark=o] table[x=dp,y=p015] {gk08konvmiddle.dat};
\addplot[color=brown,mark=triangle] table[x=dp,y=p005] {gk08konvmiddle.dat};
\legend{{$\tau=0.3\;\mathrm{s},\;h$-conv.},{$\tau=0.15\;\mathrm{s},\;h$-conv.},{$\tau=0.05\;\mathrm{s},\;h$-conv.},
{$\tau=0.3\;\mathrm{s},\;p$-conv.},{{$\tau=0.15\;\mathrm{s},\;p$-conv.}},{$\tau=0.05\;\mathrm{s},\;p$-conv.}}
\end{loglogaxis}
\end{tikzpicture}
\end{tabular}
\caption{Convergence histories of the relative errors measured in maximum norm for the front- and rear side temperature, as well as the mid-side heat flow -- GK model with $\kappa^2=0.8\;\mathrm{m}^2$}\label{fig:GK-large}
\end{figure}

\section{Discussion}
Here we presented a novel numerical approach, specifically developed for the Guyer--Krumhansl equation as that model could have significant practical interest in engineering based on the available experimental data. We successfully demonstrated that the solutions converge, thus are stable and consistent and reproduce the experimentally observed required characteristics, on contrary to COMSOL \cite{RietEtal18}. That is essential, since the usual approaches do not work, it is not possible to use the temperature as a single field variable for heat flux BCs \cite{Kov22a}.

Moreover, that technique is incredibly more efficient than the one offered by COMSOL, runs about 4 magnitudes faster thanks to the low number of DOFs and elements. That advantage becomes more significant in two and three-dimensional problems, therefore our next aim is to improve and implement the present approach for complex geometries in higher spatial dimensions. Additionally, due to the global energy crisis and chip shortage, the efficiency of algorithms has increasing importance. That huge increase in speed would enable the real-time monitoring of heterogeneous materials with complex inner structure without massive computational capacity.

\section{Acknowledgements}
The research reported in this paper is part of project no. BME-NVA-02, implemented with the support provided by the Ministry of Innovation and Technology of Hungary from the National Research, Development and Innovation Fund, financed under the TKP2021 funding scheme.
This paper was supported by the János Bolyai Research Scholarship of the Hungarian Academy of Sciences.
The research reported in this paper and carried out at BME has been supported by the grants National Research, Development and Innovation Office-NKFIH FK 134277.


\end{document}